\documentclass[12pt]{article}
\usepackage{amsmath, amssymb,amsthm}

\begin{document}

\title{\bf On mod $p$ singular  modular forms II}
\author{Siegfried B\"ocherer and Toshiyuki Kikuta}
\maketitle

\noindent 
{\bf 2020 Mathematics subject classification}: Primary 11F33 $\cdot$ Secondary 11F46\\
\noindent
{\bf Key words}: Siegel modular forms, Congruences for modular forms, Singular\\

\begin{abstract}
  We generalize the notion of mod $p^m$ singular Siegel modular forms of $p$-rank $r$
  to the vector-valued case and we show that also in this case a congruence mod $(p-1)p^{m-1}$ 
  between the scalar weight and the $p$-rank must hold.   
In some sense our proof is even simpler than the one we gave previously in the scaler valued case. 
\end{abstract}  

\section{\bf Introduction}

 Integrality properties of Fourier coefficients and of Hecke eingenvalues
and congruences among them are  important
topics in the arithmetic theory of modular forms. 
The recent impressive works of \cite{At}  show that one
should consider such questions also in the context
of vector-valued Siegel modular forms.
\\
We study the following question: \\
{\it Can it happen that a degree n Siegel modular form with integral Fourier coefficients
  has all its rank n Fourier coefficients divisible by $p$ if there is at least a rank n-1
  Fourier coefficient coprime to $p$?}\\
In \cite{BoKi} we studied this problem for the case of
scalar-valued automorphy factors.
We showed that this is only possible if $n$ and the weight $k$
are tied together by a
congruence mod $p-1$. Our aim is a generalization of this statement
for the vector-valued case.
We emphasize that our method is different from (and in some sense simpler than) the
one used in \cite{BoKi}, where the
theta decomposition of Jacobi forms (applied to Fourier-Jacobi coefficients of $f$)
was an essential tool; we avoid such
decomposition, which is somewhat delicate in the vector-valued case).\\
We adopt some techniques from the theory of
vector-valued singular modular forms over ${\mathbb C}$
as created by Freitag \cite{Frei} and from our work on (scalar-valued) mod
$p^m$ singular
modular forms \cite{BoKi2}.
\\
Here is our main result (notations will be explained below).\\[0.3cm]
{\bf Theorem:} {\it Let $\rho$ be an irreducible
  $\ell$-dimensional polynomial representation of
$GL(n,{\mathbb C})$, of scalar weight $k$
and realized integrally on ${\mathbb C}^{\ell}$ and $p$ an odd prime.
Let $f$ be an integral  modular form of automorphy factor $\rho$
for a group $\Gamma$ containing $\Gamma_0(p^s)\cap \Gamma_1(N)$
with $ N$ coprime to $p$, possibly with real nebentypus mod $p$.
If $f$ is mod $p^m$ singular of rank $r<n$ then we have
$$2k-r \equiv 0\bmod (p-1)p^{m-1}.$$
}

For simplicity, the theorem above is stated for the full modular group.  
In the proof we will cover a more general situation, allowing
certain congruence subgroups and also quadratic nebentypus characters mod
$p$. For even more general situations, we refer to the last section.

\section{Preliminaries}
\subsection{Siegel modular forms}

Let $n$ be a positive integer,
${\mathbb H}_n$ the Siegel upper half space of degree $n$ and $\Gamma^n$ the
Siegel modular group of degree $n$. For two coprime positive integers
$N$ and $N'$ we will consider
congruence subgroups of type
$$\Gamma^n_0(N,N'):=\{\left(\begin{array}{cc} A & B\\
  C & D\end{array}\right)\in \Gamma^n\mid C\equiv 0\bmod NN',
  \det(D)\equiv 1\bmod N'\}.$$
In standard notation, $\Gamma^n_0(N,N')=\Gamma^n_0(N)\cap \Gamma^n_1(N')$.\\   
  Let $\rho:GL(n,{\mathbb C})
    \longrightarrow GL(V) $ be an
    irreducible $\ell$-dimensional polynomial representation, we will always use a version in coordiantes, i.e.,  $V={\mathbb C}^{\ell}$.
    We define the slash-operator for
    $V$-valued functions on ${\mathbb H}_n$ and
    $M=\left(\begin{array}{cc} A & B\\
      C & D\end{array}\right)\in Sp(n,{\mathbb R})$ by
    $$F\mid_{\rho}M =\rho(CZ+D)^{-1}F(M<Z>).$$

    For a Dirichlet character $\chi$ mod $N$, a number $N'$ coprime to $N$
    and a representation $\rho$ as above,
we define the space of modular forms of automorphy factor $\rho$
for $\Gamma:=\Gamma^n_0(N,N')$ as
$$M^n_{\rho,\chi}(\Gamma):= \{F:{\mathbb H}_n\longrightarrow V\mid
\forall \gamma\in \Gamma: F\mid_{\rho}g=\chi(\det(D))\cdot F, \,\,F\,
\mbox{holomorphic}\}
$$
(and the standard additional condition for $n=1$). 

For general properties of Siegel modular forms we refer to the books
\cite{Frei,Kli}.

We use an ``integral version'' of $\rho$, i.e.,
we assume that - after choosing appropriate coordinates in the representation
space $V$ - we have
$$\rho(GL(n,{\mathbb Z}))\subset GL(\ell,{\mathbb Z}).$$
In particular, the polynomial functions defining
$\rho$ have coefficients in ${\mathbb Q}$.
Such a choice is always possible, see \cite{Green}.

Any modular form $F\in M^n_{\rho,\chi}(\Gamma)$ has a Fourier expansion
$$F(Z)=\sum_{
  T\in \Lambda^n} a_F(T) e^{2\pi i tr(TZ)},$$
  where
  $$\Lambda^n=\{T=(t_{ij})\in Sym_n({\mathbb Q)}\,\mid t_{ii},\ 2t_{ij}\in {\mathbb Z},\ T \,\mbox{positive\ semidefinite}\}.$$
  Furthermore, for $0\leq r\leq n$ we denote by
  $\Lambda^n_r$ the subset of all $T$ with $rank(T)=r$.
  \\
 Later on we also need a Fourier-Jacobi expansion of $F$; for that purpose we  
 we decompose $Z\in {\mathbb H}_n$ into block matrices
    $Z=\left(\begin{array}{cc} z_1 & z_2\\
   z_2^t & z_4\end{array}\right)$ with $z_4\in {\mathbb H}_r$
   and we get a partial Fourier series

   $$F(Z) = \sum_{T\in \Lambda^r}\Phi_T(z_1,z_2)e^{2\pi i tr(Tz_4)}.$$
 
Under the integrality assumption from above for $\rho$
  we can talk about integral Siegel modular forms by
requesting all Fourier
coefficients to be in ${\mathbb Z}^{\ell}$.

We also mention that $M^n_{\rho,\chi}=M^n_{\rho,\chi}({\mathbb Z})\otimes {\mathbb C}$ 
always holds provided that $\chi$ is a quadratic character (see \cite{Bovec} for a simple proof in the level one case),
but such a property will not be
needed in our work.
  \\
  We call such a modular form mod $p^m$ singular of rank $r$ if all
  Fourier coefficients
  $a_F(T)$ with $rank(T)>r$ are in $p^m\cdot {\mathbb Z}^{\ell}$
  and there exists some Fourier coefficient of rank $r$ which is
  not congruent zero mod $p$.\\
  By  general representation theory there is a weight
 $${\bf k} = (k_1\geq k_2\geq \dots \geq k_n)$$ 
  associated to $\rho$ (``highest weight'').
  \\
  A mod $p^m$ singular modular form $F$ of rank $r$ is in
  particular a noncusp form,
  it remains nonzero after $(n-r)$-fold application of
  Siegel's $\Phi$-operator; moreover, $\Phi^{n-r+1}(F)$ is a 
  mod $p^m$ singular modular form of rank $r$ and degree $r+1$.\\  
  A theorem of Weissauer \cite{Weiss} tells us that
  $$k_n=\dots = k_{n-r}$$
  must hold. We call this $k:=k_n$ the (scalar) weight of $\rho$
  and we write $\rho=\rho_0\otimes \det^k$.
  Note that $\rho_0$ is still polynomial.\\

  \section{On some partial Fourier series of $F$}
  We study some subseries of the Fourier expansion of a modular form,
  defined by the rank of $T$ and by divisibility properties of
  the entries of $T$ (respectively). This section is independent
  of congruence properties.

  \subsection{The subseries defined by fixing the rank of $T$}
  We want to avoid the use of the theta decomposition of the Fourier-Jacobi coefficients of $f$,
  which was used in \cite{BoKi} as a main tool.\\
  Instead we start, for an arbitrary modular form
  $F\in M^n_{\rho,\chi }(\Gamma)$  with a description of the
  ``rank $r$''-part of the
  Fourier expansion, following the procedure in \cite{BoKi2}:

    We express the partial Fourier series
    $$F^*:=\sum_{T\in \Lambda_n, rank(T)=r} a_F(S) e^{2\pi i tr(TZ)}$$
  by
  $$F^0:=\sum_{T\in  \Lambda^r_r} a_F(\left(\begin{array}{cc} 0 & 0\\
    0 & T\end{array}\right)) e^{2\pi i tr(Tz_4)},$$
    using the relation
    $$a_F(UTU^t)=\rho(U)\chi(\det(U))a_F(T),$$    
    which holds for all $T$ and all $U\in GL(n,{\mathbb Z})$ with $\det(U)\equiv 1\bmod N'$.
    \\
    We obtain, using the submatrices $z_1,z_2,z_4$ introduced above,
   $$F^*= \sum_U \sum_{T\in \Lambda^r_r} \rho(U) a_F(\left(\begin{array}{cc} 0 & 0\\
      0 & T\end{array}\right)) e^{2\pi i tr( U\cdot \left(\begin{array}{cc}
          0 & 0\\
          0 & T\end{array}\right)U^t\cdot Z)}.$$
    
      The expression in the exponential is then equal to

      $$tr(u_2Tu_2^{t}z_1+2 u_2Tu_4z_2) +tr(u_4Tu_4^{t}z_4).$$
    
      The summation over $U$ goes over $GL(n,{\mathbb Z})/ GL(n,r,{\mathbb Z})$
      with $GL(n,r,{\mathbb Z})$ being defined by

 $$ \{U=\left(\begin{array}{cc} u_1 & 0\\
        u_3 & u_4\end{array}\right)\,\mid
\begin{array}{l}
        u_1\in GL(n-r,{\mathbb Z}), u_3\in
        {\mathbb Z}^{r,n-r}, u_4\in Aut(T),\\
        \det(u_1)\cdot\det(u_4)\equiv 1\bmod N'\end{array}
        \}.$$

      \subsection{Imposing congruence conditions on the upper right block submatrix of $T$}
      In the context of the doubling method we already used a kind of ``partial twists'' of modular forms
      to modify the primitive nebentypus character of the modular form in question, see \cite{BoSchm}.
      In loc.cit., remark 2.2. we mentioned that this also works
      for trivial character.
      We will use it in the following context: \\
   
      \noindent   
 {\bf Observation:} {\it Let $F$ be a modular form of
   degree n (possibly vector-valued and with nebentypus character mod $N$)
   of level
    $\Gamma_0(N)\cap \Gamma_1(N')$ with Fourier expansion
  $$F(Z)= \sum_S a(S)  e^{2\pi itr(TZ)}.$$
  Then for any $t\geq 0$ the subseries
  $\tilde{F}(Z)$, where $S$ runs only over those matrices, which are of the form
  $$S=\left(\begin{array}{cc} S_1 & S_2\\
    S_2^t & S_4\end{array}\right) \quad\mbox{with}\quad  S_1\in \Lambda_{n-r},\ S_4\in \Lambda_r
    \quad \mbox{and} \quad S_2\equiv 0\bmod p^t$$
    is still modular, in particular for
    $$\Gamma^{n-r}_0(N\cdot p^{2t})\cap \Gamma_1(N')\times
    \Gamma^r_0(N\cdot p^{2t})\cap \Gamma_1(N')
    \subset \Gamma_0^n(N)\cap \Gamma_1(N')$$
Also the nebentypus character $\chi$ is maintained. 
  }\\[0.3cm]
  Here we identify the product $Sp(n-r)\times Sp(r)$ with a subgroup of $Sp(n)$
via 
the natural embedding defined by
$$\left(\begin{array}{cc} A & B \\C & D\end{array}\right)\times
  \left(\begin{array}{cc} a & b\\ c & d\end{array}\right) \hookrightarrow
  \left(\begin{array}{cc}
    \begin{array}{cc} A & 0\\
      0 & a
    \end{array} &\begin{array}{cc} B & 0\\
      0 & b
    \end{array}\\
    \begin{array}{cc} C & 0 \\
      0 & c
    \end{array} & \begin{array}{cc} D & 0\\
      0 & d
  \end{array}\end{array}\right).
  $$

  This follows from a standard matrix calculation taking into account
  that going from $F$ to $\tilde{F}$ is achieved by

  $$F\longmapsto \sum_R F\mid_{\rho} \left(\begin{array}{cc} 1_n &
    \begin{array}{cc} 0 & \frac{R}{ p^t}\\
      \frac{R^t}{ p^t} & 0 \end{array}\\[0.4cm]
    0_n & 1_n\end{array}\right).$$
      
The summation goes over all elements of ${\mathbb Z}^{(n-r,r)}$ mod $p^t$.

 \section{Proof of theorem}     
   
      We now apply the consideration from above to a modular form
      $F\in M^n_{\rho,\chi}(\Gamma)_{\mathbb Z}$, which is mod $p^m$ singular
      of rank $r$.
We extract for this $f$ the Fourier-Jacobi coefficient mod $p^m$
for $T\in \Lambda^r_r$ with
$a(\left(\begin{array}{cc} 0 & 0\\
  0 & T\end{array}\right) ) $ not congruent $0$ mod $p$ and $\det(T)$
  minimal with this property.
Then the summation over $U$ just becomes summation over $u_2$.

For such $T$ we get

\begin{equation}
\Phi_T(z_1,z_2)\equiv \sum_{u_2} \rho(\left(\begin{array}{cc}
  1 & u_2\\
  0 &1\end{array}\right))\cdot  a(\left(\begin{array}{cc} 0 & 0\\ 0 & T
  \end{array}\right) ) e^{2\pi i tr(u_2Tu_2^t z_1 + 2 u_2Tz_2)}\bmod p^m.
  \label{theta} \end{equation}

Everything would be much easier if we could restrict the summation to
those $u_2$ for which $\rho(\left(\begin{array}{cc} 1 & u_2\\
  0 & 1\end{array}\right))$
  would just be congruent mod $p^m$ to the identity.
  
  To achieve this, we use the observation from the previous section to
  switch from $F$ to $\tilde{F}$ with sufficiently large $t$.
     This is still singular mod $p^m$ of rank $r$.

     The condition
``$S_2\equiv 0$ mod $p^t$'' from the observation now becomes 
     ``$T\cdot u_2\equiv 0$ mod $p^t$'',
     which 
can be rephrased as
$$T\cdot u_2\in T\cdot {\mathbb Z}^r\cap p^t\cdot {\mathbb Z}^r$$
i.e., 
\begin{equation} u_2\in {\mathbb Z}^r\cap p^t\cdot T^{-1} \cdot {\mathbb Z}^r= R\cdot {\mathbb Z}^t\label{sublattice}\end{equation}
  for a suitable integral $r\times r$ matrix $R$ of maximal rank.
  Then for sufficiently large t (depending on $T$)
  we have $u_2\equiv 0$ mod $p^m$ for such $u_2$.
  \\
  Now we look at the Fourier-Jacobi coefficient
  $\tilde{\Phi}$ mod $p^m$ for $\tilde{F}$.
  
  The expression in (\ref{theta}) for
  $\tilde{\Phi}_T(z_1,0)$ then runs only over $u_2$ 
with $u_2\equiv 0$ mod $p^t$; on such $u_2$ the value of $\rho$ is just
the identity mod $p^m$.

Then we obtain
$$\tilde{\Phi}_T(z_1,0)\equiv \sum_{u_2\in R\cdot {\mathbb Z}^r} a(\left(\begin{array}{cc} 0 & 0\\
  0 & T\end{array}\right) e^{2\pi i tr(u_2Tu_2^tz_1)}\bmod p^m$$
  We may rewrite this, using (\ref{sublattice}) and putting ${\mathcal R}:= R^tTS$
  
\begin{equation}
\tilde{\Phi}_T(z_1,0)\equiv \theta^{(n-r)}_{\mathcal R}(z_1)\cdot a(\left(\begin{array}{cc} 0 & 0\\
    0 & T\end{array}\right))\bmod p^m\label{congruence1}\end{equation}
with the theta series
$$\theta_{\mathcal R}^{(n-r)}(z_1):=\sum_{X\in {\mathbb Z}^{(r,n-r)}}
e^{2\pi itr(X^t{\mathcal R}X\cdot z_1)}.$$

The change from $T$ to ${\mathcal R}:=S^t\cdot T\cdot S$ comes from
the condition (\ref{sublattice}).

  As in the scalar-valued case we may assume $n=r+1$; this is not really necessary but the ``branching laws'' for $\rho$ are somewhat easier to handle.
  \\
  To understand $F|_{z_2=0}$ and $\Phi_T(z_1,0)$ as modular forms,
  we decompose the representation space $V$ for $\rho_0$ as a $GL(1)\times GL(n-1)$-module. The branching laws (e.g., Goodman-Wallach
  \cite[Theorem 8.1.2]{GoWa}
  )
  give
  $$V=\oplus V_{i}\otimes W_{i}$$
  where $V_j$ is one dimensional with $x\in GL(1,{\mathbb C})$ acting on
  $V_{i}$ by multiplication by $x^{i}$.
  We do not need the explicit description of the $GL(n-1,{\mathbb C})$-module $W_{i}$ here.
  
  Using the elementary divisor theorem, we may choose
  a basis
  $${\mathfrak a}_1,\dots , {\mathfrak a}_m$$
  and natural numbers $\alpha_1\mid \dots \mid \alpha_r$
  such that
  ${\mathfrak a}_1,\dots , {\mathfrak a}_m$ is a ${\mathbb Z}$-basis
  of $V({\mathbb Z})$ and
  $\alpha_1\cdot {\mathfrak a}_1,\dots , \alpha_r\cdot{\mathfrak a}_r$
  is a ${\mathbb Z}$-basis of $\left(V_0\otimes W_0\right)\cap V({\mathbb Z})$.
  We do not claim that $r$ is the dimension of the ${\mathbb C}$-vector space
  $V_0\otimes W_0$.
\\[0.3cm]
  How can we interpret the congruence (\ref{congruence1})?
  
  We observe that $\tilde{\Phi}_T(z_1,0)$ is a $V$-valued function,
  more precisely
  $$
  \tilde{\Phi}_T(z_1,0)=\sum_i g_i(z_1)\cdot {\mathfrak a}_i$$
  where the $g_i$ are modular forms of weight $k$ for
  $i\leq r$ and of higher weight otherwise (all of degree one).

  On the other hand,
   $a(\left(\begin{array}{cc} 0 & 0\\
      0 & T\end{array}\right) )$ is $(V_0\otimes W_0)\cap V({\mathbb Z})$-valued and therefore
  $$\theta_{\mathcal R}(z_1)\cdot a(\left(\begin{array}{cc} 0 & 0\\
        0 & T\end{array}\right) )=
        \theta_{\mathcal R}(z_1)\sum_{j=1}^r \beta_j\cdot \alpha_j{\mathfrak a}_j$$
      
        with integer coefficients $\beta_j$ and at least
        one $\beta_{j_0}\cdot\alpha_{j_0}$ not congruent zero mod $p$.

  The congruence (\ref{congruence1}) then implies that for such $j_0$ the congruence
  $$ g_{j_0}\equiv \beta_{j_0}\cdot \alpha_{j_0}\cdot \theta_{\mathcal R}\bmod p$$
  holds.
  \\
  This is a congruence mod $p^m$ between a degree one modular form $g:=g_{j_0}$
  of level $\Gamma_1(N)\cap \Gamma_0(p^{2t})$ with possible quadratic nebentypus
  mod $p$ and a theta series $\theta_{\mathcal R}$ of weight $\frac{r}{2}$.
  To avoid the case of half-integral weights, we consider the congruence
  between $g^2$ and the theta series for ${\mathcal R}\perp {\mathcal R}$.
  This theta series of weight $r$ then may have
  nebentypus $\left(\frac{(-1)^r\det(R)^2}{.}\right) $ which is a nontrivial
  character at most modulo $4$.
  After possibly applying level change mod $p^m$ to $g^2$ and
  the $\theta$-series (their weights would not change mod $(p-1)p^{m-1}$ ,
  we arrive at a congruence between two modular forms
  of level coprime to $p$. We may now argue as in
  \cite[Corollary 4.4.2]{Katz}, see also
  \cite[section 5]{BoKi}.
  We obtain the desired congruence
  $$2k\equiv r\bmod (p-1)p^{m-1}.$$

\section{Final remarks:}

\begin{itemize}
    
\item In this note (mainly for simplicity of notation) we
  stick to congruences for
powers of rational primes $p\not=2$, but in the same style as in \cite{BoKi}
it is possible to consider a situation
over arbitrary number fields and to
consider congruences modulo powers of prime ideals,
including the dyadic case.
Also half-integral weights can be considered.
As for nebentypus characters, arbitrary characters with conductors coprime
to $p$ are (implicitly) allowed by our choice of $\Gamma=\Gamma_0(p^s)\cap \Gamma_1(N')$,
but we  only allow quadratic characters
in its  $p$-part. We refer to \cite{BoKi} for details.\\
In addition, we just mention that already for scalar-valued
forms congruences for modular forms
of nonquadratic characters modulo arbitrary powers of $p$ (or mod powers of
prime
ideals over $p$ are rather
delicate, in particular for mod $p$-power singular forms,
see \cite{BoKi3}.
\item Irreducibility of $\rho$ is not really needed; indeed the proof from
  above works (for $n=r+1$) in the case of arbitrary $\rho$, provided that
  all the irreducible components are of the same scalar weight $k$.
  In more general cases we would have to request that the decomposition of $\rho$ into irreducible components can be ``defined over ${\mathbb Z}$''.
 \item  There might also be a possibility of generalizing our
  structure theorem on mod $p^m$ singular forms \cite{BoKi2}
  to the vector-valued case. For that purpose however we cannot just
  switch to ``$T\cdot u_2\equiv 0$ mod $p$'' as
  above; we would have to work directly
  with series of type (\ref{theta}).
\end{itemize}

\section*{Acknowledgment}
This work was supported by JSPS KAKENHI Grant Number 22K03259.

  Siegfried B\"ocherer\\
  Kunzenhof 4B\\
  79117 Freiburg, Germany\\
  E-mail: boecherer@t-online.de\\[0.2cm]
  Toshiyuki Kikuta\\
  Department of Information and Systems Engineering\\
  Faculty of Information Engineering\\
  Fukuoka  Institute of Technology\\
  3-30-1 Wajiro-higashi, Higashi-ku, Fukuoka 811-0295, Japan\\
  E-mail: kikuta@fit.ac.jp
  
\end{document}